\documentclass{amsart}
\usepackage{amssymb,mathrsfs,amsmath,amsthm,color,bm,mathtools,bbm,wasysym,cases,mathdots}
\usepackage[pagebackref]{hyperref}
\usepackage[noadjust]{cite}

\usepackage{enumerate}
\usepackage[centering]{geometry}
\newtheorem{theorem}{Theorem}[section]
\newtheorem{lemma}{Lemma}[section]
\newtheorem{proposition}{Proposition}[section]

\theoremstyle{definition}
\newtheorem{definition}{Definition}
\newtheorem{remark}{Remark}

\numberwithin{equation}{section}

\newcommand{\R}{\mathbb R}
\newcommand{\Q}{\mathbb Q}
\newcommand{\Z}{\mathbb Z}
\newcommand{\N}{\mathbb N}
\newcommand{\T}{\mathbb T}
\newcommand{\Mat}{\mathrm{M}}
\newcommand{\Exact}{\operatorname{Exact}}
\newcommand{\supp}{\operatorname{supp}}

\newcommand{\dif}{{\,\mathrm d}}
\newcommand{\ve}{\varepsilon}

\newcommand{\cq}{\mathcal Q}

\newcommand{\cd}{\mathcal D}
\newcommand{\ce}{\mathcal E}

\newcommand{\cL}{\mathcal L}

\newcommand{\bb}{\mathbf b}
\newcommand{\bq}{\mathbf q}
\newcommand{\bp}{\mathbf p}
\newcommand{\br}{\mathbf r}
\newcommand{\bu}{\mathbf u}

\newcommand{\bt}{\mathbf t}
\newcommand{\bk}{\mathbf k}
\newcommand{\bzero}{\mathbf 0}

\newcommand{\bta}{\bm\theta}

\newcommand{\fdim}{\dim_{\rm F}}
\newcommand{\hdim}{\dim_{\rm H}}

\begin{document}

\title[Fourier dimension of exact approximation]
{Fourier dimension for exactly approximable matrices}

\author[Y. He]{Yubin He}
\address[Y. He]{School of Mathematical Sciences, Shenzhen University, 518060, Shenzhen, China}
\email{ybhe@szu.edu.cn}

\author[B. Li]{Bixuan Li}
\address[B. Li]{Institute of Mathematics, Academy of Mathematics and System Science, CAS, Zhongguancun East Road 55, Beijing 100190, P. R. China}
\email{libixuan@amss.ac.cn}

\author[B. Wang]{Baowei Wang}
\address[B. Wang]{School of Mathematics and Statistics, Huazhong University of Science and Technology, Wuhan 430074, P. R. China}
\email{bwei\_wang@hust.edu.cn}


\begin{abstract}
We determine the Fourier dimension of exactly approximation set for all dimensions. More precisely, let $m,n\ge1$, $\bta\in\R^m$, and $\psi:\N\to \R^+$ be a non-increasing approximate function. Denote by $\Exact_{m,n}(\psi,\bta)$ the set of $m\times n$ matrices which are
inhomogeneously $\psi$-approximable with shift $\bta$, but not
$c\psi$-approximable for any fixed $0<c<1$.
We prove that
$$\dim_{\text{F}}
\Exact_{m,n}(\psi,\bta)=\frac{2n}{1+\tau_\psi},\ \text{whenever }\tau_{\psi}>\frac{n}{m},$$
where $\tau_{\psi}$ is the lower order of $\psi$ at infinity.
\end{abstract}

\maketitle

\section{Introduction}
For a positive function $\psi:\mathbb N\to \R^+$, the set of $\psi$-well-approximable numbers is defined by
\[
W(\psi)=\left\{x\in[0,1]: |qx-p|\leq\psi(q)\text{ for infinitely many }(p,q)\in\mathbb Z\times\mathbb N\right\}.
\]
Khintchine's theorem \cite{Khintchine} and Jarn\'{i}k's theorem \cite{Jar1} on Lebesgue measure and Hausdorff dimension of $W(\psi)$ were the starting point of metric Diophantine approximation. These classical results were subsequently extended to higher-dimensional settings, including simultaneous approximation and, more generally, systems of linear forms.

In this paper, we work on the setting of system of linear forms, so fix some notation at first. For a vector $\mathbf x=(x_1,\ldots,x_d)^\top\in\mathbb R^d$, let
\[
|\mathbf x|=\max_{1\le i\le d}|x_i|,\qquad \|\mathbf x\|=\min_{\mathbf k\in\mathbb Z^d}|\mathbf x-\mathbf k|
\]
 denote its supremum norm and its distance to the integers $\mathbb Z^d$, respectively. Let $\Mat_{m,n}(\mathbb R)$ be the spaces of $m\times n$ matrices with entries in $\mathbb R$.

Throughout this paper, let $\boldsymbol\theta\in\mathbb R^m$ and $\psi:\mathbb N\to \R^+$ be a non-increasing function. The set of inhomogeneous $\psi$-well approximable matrices is defined by
$$W_{m,n}(\psi,\boldsymbol\theta)=\left\{A\in\Mat_{m,n}(\mathbb T):\|A\mathbf q-\boldsymbol\theta\|\leq\psi(|\mathbf q|)\text{ for infinitely many }
\mathbf q\in\mathbb Z^n\setminus\{\mathbf 0\}
\right\}.$$

The metric theory for $W_{m,n}(\psi,\boldsymbol\theta)$ has been well established in the literature. The measure theory for the homogeneous case ($\boldsymbol\theta=\mathbf 0$) was presented by Khintchine and Groshev \cite{Gro}, and the inhomogeneous case was given by Sprind\v{z}uk \cite{Sprindzuk}. Hausdorff dimension were obtained by Dodson \cite{Dodson} for the homogeneous case and by Levesley \cite{Levesley} for the inhomogeneous case by showing
$$\dim_{\mathrm H}W_{m,n}(\psi,\boldsymbol\theta)=m(n-1)+\frac{m+n}{1+\tau_\psi},\ \text{where } \tau_\psi=\liminf_{q\to\infty}\frac{-\log\psi(q)}{\log q}.$$
For further related results, one is referred to \cite{BerDV,Ber,Bug,DicV} and the references therein.

In classic analysis, a widely concerned question is whether a measure zero set $E$ can support a Rajchman measure, i.e. $\exists \ \mu$ supported on $E$ such that $$
\widehat\mu(\xi):=\int e^{2\pi i x\cdot\xi}d \mu(x)\to 0, \ \ {\text{as}}\ \ |\xi|\to \infty.
$$ It can be asked further about the maximal Fourier decay of measures supported on $E$. In fact, it also has connection with the
Hausdorff dimension of $E$ through the Fourier representation of Riesz energies (see \cite[Theorem~3.10]{Mattila}). For \(0<t<d\), the \(t\)-energy of a finite positive Borel measure \(\mu\) on \(\mathbb R^d\) satisfies
$$
I_t(\mu)
:=\iint_{\mathbb R^d\times\mathbb R^d}
\frac{\dif\mu(x)\,\dif\mu(y)}{|x-y|_2^t}
=c_{d,t}\int_{\mathbb R^d}
|\widehat\mu(\xi)|^2|\xi|_2^{t-d}\,\dif\xi,
$$
where \(|\cdot|_2\) denotes the Euclidean norm and \(c_{d,t}>0\) depends only on \(d\) and \(t\). In particular, if
$$
|\widehat\mu(\xi)|\ll (1+|\xi|_2)^{-s/2}\asymp(1+|\xi|)^{-s/2},
$$
 then \(I_t(\mu)<\infty\) for every \(0<t<s\). Consequently, by \cite[Theorem~4.13]{Falconer}, for any Borel set \(E\), $$
\hdim E\ge s, \ {\text{if}}\ \mu(E)>0.
$$ 

 This observation motivates the notion of Fourier dimension, which quantifies the maximal Fourier decay attainable by measures carried by a set.

\begin{definition}[Fourier dimension]
Let $E\subset\R^d$ be a Borel set. The \textit{Fourier dimension} of $E$ is defined by
\begin{equation*}\label{eq:fourier-dimension}
\dim_{\text{F}} E=\sup\Bigl\{s\in[0,d]:\ \text{$\exists\ \mu\in\mathcal{P}(E)$ such that }|\widehat\mu(\xi)|
\ll(1+|\xi|)^{-s/2}\text{ for all }\xi\in\R^d\Bigr\},
\end{equation*}
where $\mathcal{P}(E)$ denote the collection of Borel probability measures on $\mathbb{R}^d$ that give full measure to $E$.
\end{definition}

\begin{remark}
Different notions of Fourier dimension are used in the literature, often without an explicit distinction. Here, we follow the definition given by Ekstr\"om, Persson and Schmeling \cite{EPS}. Another widely used definition is called compact Fourier dimension \(\dim_{\mathrm{FC}}E\) in \cite{EPS}, which is given as
$$
\dim_{\mathrm{FC}}E=\sup\bigl\{\dim_{\mathrm F}K:K\subset E\text{ is compact}\bigr\}.
$$

There are some differences between these two definitions (see \cite{EPS}). However, there is no reason to say which one would be better. As far as the metric theory is concerned, they almost share the same status: \begin{itemize}
  \item both of them provide the lower bound of the Hausdorff dimension $$
  \dim_{\mathrm{FC}}E\le \hdim E, \ \ \dim_{\mathrm F}E\le \hdim E.
  $$

  \item positiveness of them guarantees the existence of normal numbers in $E$ by applying Davenport--Erd\"{o}s--LeVeque criterion \cite{DEL}.

  \item similar metric result when studies the Diophantine approximation on $E$, for example $$
  \{x\in E: \|qx-\theta\|<\psi(q), \ {\text{for infinitely many}}\ n\in \N\}
  $$ with respect to the measure supporting on $E$ or giving full measure to $E$ with polynomial Fourier decay \cite{PollV,TanZhou3}.
\end{itemize}
The flexibility of not requiring the measure fully supported on $E$ gives us more freedom on the choice of the measure and enables us to deal
with more complex setting, for example exact approximation as will be addressed later.
\end{remark}

The study of the Fourier dimension was originated from a question asked by Beurling: whether there exists set whose Fourier dimension equals to its Hausdorff dimension. The first example (a random set) was constructed by Salem \cite{Sam}, so a set with such a property is called {\em Salem set} now. The first deterministic Salem set was found by Kaufman \cite{Kaufman81} in Diophantine approximation, where he showed that the set of well-approximable numbers is a Salem set, i.e. for $\psi(q)=q^{-\tau}$ with $\tau>1$,
$$\dim_{\text{FC}} W_{1,1}(\psi,0)
=\dim_{\text{H}} W_{1,1}(\psi,0)=\frac{2}{1+\tau}.$$

Building on Kaufman's work, Bluhm \cite{Bluhm98} developed a deterministic Cantor-type construction of Salem sets. Hambrook \cite{Ham1} subsequently extended this approach to general approximation functions, without imposing monotonicity assumption, and to inhomogeneous setting.

In higher dimensions $mn>1$, the Salem property does not persist in general. Hambrook and Yu \cite{HamY} determined the compact Fourier dimension
of well-approximable matrix sets, 
which shows, except in the one-dimensional case, these sets are non-Salem. For inhomogeneous setting and approximation functions without monotonicity assumption, Cai and Hambrook \cite{CaiH} showed that
$$
\fdim W_{m,n}(\psi,\boldsymbol\theta)
=\dim_{\mathrm{FC}}W_{m,n}(\psi,\boldsymbol\theta)=\min\{2s_\psi,mn\},  \ \ {\text{when}}\  \sum_{q=1}^{\infty}q^{n-1}\psi(q)^m<\infty,
$$
where
$$
s_\psi=\inf\left\{s>0: \sum_{q=1}^{\infty}q^{n-1}\left(\frac{\psi(q)}{q}\right)^s<\infty\right\}.
$$
Although the first equality was not explicitly stated in their paper, it follows from their arguments. Consequently, when \(\psi\) is
non-increasing and $\tau_\psi>n/m$, one has
\begin{equation}\label{eq:fdim}
\fdim W_{m,n}(\psi,\boldsymbol\theta)=\dim_{\mathrm{FC}}W_{m,n}(\psi,\boldsymbol\theta)=\frac{2n}{1+\tau_\psi}.
\end{equation}

Further results on the compact Fourier dimension in weighted and multiplicative settings, see \cite{He,TanZhou1} and inhomogeneous
Duffin--Schaeffer-type sets with coprimality constraints, see \cite{TanZhou2}.
See also explicit Salem sets in $\R^d$ for $d\ge 2$ in \cite{FraH, Ham}.

In this paper, we focus on the exact approximation set introduced by Jarn\'{i}k \cite{Jarnik}. Recall that a matrix $A\in\Mat_{m,n}(\R)$ is called exactly $(\psi, \boldsymbol\theta)$-approximable, if 	
\begin{itemize}
\item there are infinitely many $(\mathbf p,\mathbf q)\in\Z^m\times(\Z^n\setminus\{\mathbf0\})$satisfying
$$|A\mathbf q-\mathbf p-\boldsymbol\theta|\le\psi(|\mathbf q|);$$

\item for every $0<c<1$, there are only finitely many pairs $(\mathbf p,\mathbf q)\in\Z^m\times(\Z^n\setminus\{\mathbf0\})$ satisfying
$$|A\mathbf q-\mathbf p-\boldsymbol\theta|\le c\,\psi(|\mathbf q|).$$
	\end{itemize}
Denote the set of exactly $(\psi, \boldsymbol\theta)$-approximable matrices by
$$\Exact_{m,n}(\psi,\boldsymbol\theta):=W_{m,n}(\psi,\boldsymbol\theta)\setminus\bigcup_{0<c<1}W_{m,n}(c\psi,\boldsymbol\theta).$$
Because of the integer-translation invariance of $\Exact_{m,n}(\psi,\boldsymbol\theta)$, we focus on the $m\times n$ matrices $\Mat_{m,n}(\mathbb T)$ with entries in  $\mathbb T=\mathbb R/\mathbb Z$. We use $\mathcal L_d$ to denote the normalized Haar measure on $\mathbb T^d$.

The Hausdorff dimension of exact approximation sets has been extensively investigated. In the one dimensional case, it was studied by Bugeaud \cite{Bug1,Bu2} and Bugeaud and Moreira \cite{BugM}. Their work was extended to simultaneous approximation by Bandi and de Saxc\'e \cite{BanS}, yielding
\[
\dim_{\mathrm H}\Exact_{m,1}(\psi,\mathbf0)=\dim_{\mathrm H}W_{m,1}(\psi,\mathbf0)=\frac{m+1}{1+\tau_\psi}
\]
whenever $\tau_\psi>1/m$. Recently, He, Li and Wang \cite{HLW} established the corresponding result for system of linear forms by showing
for $\boldsymbol\theta\in \R^d$, \begin{equation}\label{f1}
\dim_{\mathrm H}\Exact_{m,n}(\psi,\boldsymbol\theta)=\dim_{\mathrm H}W_{m,n}(\psi,\boldsymbol\theta)=m(n-1)+\frac{m+n}{1+\tau_\psi}
\end{equation}
whenever $\tau_\psi>n/m$. Related developments in Hausdorff measure theory and weighted exact approximation can be found in \cite{BanGN,Fregoli,TanTianZhou} and \cite{BanF}, respectively.

The (compact) Fourier dimension of exact approximation set was investigated by Fraser and Wheeler \cite{FraW,FraW1}. In \cite{FraW}, they adapted Kaufman's construction to obtain explicit positive lower bounds for large approximation exponents, together with an inhomogeneous version under Diophantine conditions on the shift. Their subsequent work \cite{FraW1} treated smaller approximation exponents by combining a continued fraction construction with oscillatory estimates. These two results establish positive (compact) Fourier dimension for $\Exact_{1,1}(q^{-\sigma},0)$ in the respective ranges
\[
\tau>\frac{1+\sqrt5}{2}
\quad\text{and}\quad
1<\tau<\frac{5+\sqrt{73}}8.
\]
Since these ranges overlap, they yield
\[
\dim_{\text{FC}} \Exact_{1,1}(q^{-\tau},0)>0
\quad\text{for every }\tau>1.
\]

Despite these advances, the exact value of the (compact) Fourier dimension of exact approximation sets remains unknown, even in the one-dimensional homogeneous case.
We answer it in this paper for all \(m,n\ge 1\) and arbitrary shifts \(\boldsymbol\theta\in\mathbb R^m\).

\begin{theorem}
Let $m,n\ge1$, $\bta\in\R^m$, and let $\psi:\N\to(0,\infty)$ be non-increasing, we have
\begin{equation*}
\dim_{{\rm{F}}}\Exact_{m,n}(\psi,\bta)=\dim_{{\rm{F}}}W_{m,n}(\psi,\bta)=\frac{2n}{1+\tau_\psi},  \ \ \text{whenever}\ \tau_{\psi}>\frac{n}{m}.
\end{equation*}
\end{theorem}

\begin{remark} The main difficult to study the structure of the exactly approximable sets is that the points in $\Exact_{m,n}(\psi,\bta)$ have
to be isolated from {\em all} but finitely many the resonant sets: $$
{\text{when}} \ n>1: {\text{subspace}}\ \  A\bold q-\bold p-\bta=0; \ \ \ \ {\text{when}} \ n=1: {\text{shifted rationals}}\ \frac{\bold p+\bta}{q}.
$$ This will require information on the separation properties of these resonant sets. The Diophantine condition on the shift $\bta$ asked by
Fraser and Wheeler \cite{FraW} is used to guarantee this.

As said before, the definition of Fourier dimension we adopted here gives us some freedom: it
is sufficient if we are clear about the separation properties for most of the resonant sets instead of all of them.
\end{remark}

\begin{remark}\label{r1}
We give some explanation about the definition of Fourier dimension we choose. A key estimate in obtaining Fourier decay is the cancellation of exponential sums. For example
\begin{equation}\label{eq:exp}
    \sum_{p=0}^{q-1}e^{2\pi i ps/q}=0  \ \ \text{whenever}\ q\nmid s.
\end{equation}
As far as the structure of exactly approximable set $\Exact_{m,n}(\psi,\bta)$ is concerned, we have to discard some neighborhood of all resonant sets, for example, $$B\left(\frac{r+\theta}{s}, \eta\frac{\psi(s)}{s}\right), \ 0<\eta<1, \ \ \ {\text{when}}\ n=1,\ m=1.$$
If there exists some integers $p,q$ such that $$
B\left(\frac{p+\theta}{q}, \frac{\psi(q)}{q}\right)\subset B\left(\frac{r+\theta}{s}, \eta\frac{\psi(s)}{s}\right),
$$ the left ball has to be deleted in constructing subsets of $\Exact_{m,n}(\psi,\bta)$. Then some integer $p$ will be absent in the sum (\ref{eq:exp}), which will break the cancellation of exponential sums and brings the difficulty in estimating the Fourier decay of measures supported on $\Exact_{m,n}(\psi,\bta)$.

However the Hausdorff dimension of $\Exact_{m,n}(\psi,\bta)$ given in (\ref{f1}) has indicated that there are still large room left in $W_{m,n}(\psi,\bta)$ after eliminating $W_{m,n}(\eta\psi, \bta)$. So a measure supported on $W_{m,n}(\psi, \bta)$ can still give large mass to $\Exact_{m,n}(\psi,\bta)$. This is the main issue to be verified which is observed from our recent work on an average separation property of resonant sets \cite{HLW}. In such a sense, one can believe that $W_{m,n}(\psi, \bta)$ and $\Exact_{m,n}(\psi,\bta)$ will have the same Fourier dimension under the definition we choose.
\end{remark}

\subsection{Organization} In Section~\ref{s:pre}, we introduce the necessary notation and recall some estimates from \cite{HLW}. In Section~\ref{s:auxiliary}, we give the main ingredients for our measure construction, including a key estimate (Proposition \ref{p31}) controlling the measure of unions of dangerous approximation sets. Finally, in Section~\ref{s:con}, we complete the construction of the measure $\mu$, establish its Fourier decay, and show that the exact approximation set has full $\mu$-measure.

\section{Preliminaries}\label{s:pre}

We begin by introducing some notation and terminology used throughout this paper. We use $a\ll b$, $a\gg b$ and $a\asymp b$ respectively to mean that $0<a/b\le e_1$, $a/b \ge e_2>0$ and $e_2\le a/b\le e_1$ for unspecified positive constants $e_1,e_2$.
\subsection{Notation}
Recall that, by the definition of $\tau_\psi$, for any $\ve>0$ and all
sufficiently large $t$,
$$\psi(t)\le t^{-(\tau_\psi-\ve)}.$$
We fix $\ve>0$ sufficiently small so that $m(\tau_\psi-\ve)>n.$ Consequently,
\begin{equation*}
	\sum_{h\ge 1}
	h^{n-1}\psi(h)^m<\infty.
\end{equation*}

We define the exact set and the dangerous set. Let $0<\eta<1$. 
\begin{itemize}
\item For the case $n=1$. For any $q\ge1$ and $\bp\in\Z^m$, denote by
\[
R_{\bp,q}
=
\frac{\bp+\bta}{q}
\]
the \textit{resonant point} and define the associated \emph{exact cube}
\[
A_{\eta}(\bp,q)=R_{\bp,q}+\left[\eta\frac{\psi(q)}{q},\,\eta\frac{\psi(q)}{q}+\ell_M(\eta)\right]^m, \ \ {\text{where}}\ \ell_M(\eta)=(1-\eta)\frac{\psi(M)}{M}.
\]
For any $s\ge1$, $\br\in\Z^m$, define the corresponding \emph{dangerous cube}
\[
D_\eta(\br,s)
= B\bigg(R_{\br,s},
\eta\frac{\psi(s)}{s}
\bigg).
\]
\item For all the case $n\ge 1$. Define
\[
\cq_M
:=
\left\{
(\bu,q)\in\Z^n:
0<|\bu|\le\frac{M}{4},\
\frac{M}{2}\le q\le M,\
q\text{ is prime}
\right\}.
\]
For any $\bq=(\bu,q)\in\cq_M$, it is clear that $|\bq|=q$ and that $\bq$ is primitive. 
When $n=1$, the vector $\bu$ is absent, we adopt the convention $\bq=q$.

\item For all the case $n\ge 1$. For any $\bq=(\bu,q)\in\cq_M$, $\bb\in\Z^n\setminus\{\bzero\}$, and
$0<\eta<1$, define
$$\ce_\eta(\bq)
	:=
	\Big\{
	A\in\Mat_{m,n}(\T):
	 A\bq-\bta
	\in\left[
	\eta\psi(q),\eta\psi(q)+q\ell_M(\eta)
	\right]^m\ (\text{mod }\mathbb{Z}^m)
	\Big\},$$
$$\cd_\eta(\bb):=\Big\{
	A\in\Mat_{m,n}(\T):
	\|A\bb-\bta\|
	<
	\eta\psi(|\bb|)
	\Big\}.$$
When $n=1$, we adopt the convention $\bb=s$. 
It is easy to see that
\begin{equation}\label{measure}
\mathcal L_{mn}\bigl(\ce_\eta(\bq)\bigr)=q^m\ell_M(\eta)^m\asymp_m(1-\eta)^m\psi(M)^m, \ \  \ \mathcal L_{mn}\bigl(\cd_\eta(\bb)\bigr)=\bigl(2\eta\psi(|\bb|)\bigr)^m.
\end{equation}
\end{itemize}

For any $\bq=(\bu,q)\in\cq_M$, since $q\le M$ and $\psi$ is non-increasing, for any $A\in\mathcal{E}_{\eta}(\bq)$, one has
$$\eta \psi(|\bq|)\le \Vert A\bq-\bta\Vert<\psi(q)= \psi(|\bq|).$$

\subsection{Regular scales}
The following regular-scale lemma \cite[Lemma 2.1]{HLW} shows the behavior of a non-increasing function.

\begin{lemma}[Regular scale]\label{l1} Let $\psi$ be a non-increasing function. For any $\ve>0,$ there are arbitrarily large integers $M$ (called regular scale) such that \begin{equation*}
\psi(M)\ge M^{-(\tau_\psi+\ve)},\quad \psi(t)\le\psi(M)\left(\frac Mt\right)^{\tau_\psi-\ve}\text{for all }t\ge M.
\end{equation*}
\end{lemma}

We also need a general Denjoy--Koksma inequality \cite[Lemma 2.3]{HLW}. Let $\{q_j\}_{j\ge 0}$ be the denominators of the convergents of $\theta$ in its continued fraction expansion. For an integer $K\ge 1$, define
\begin{equation*}
Q_\theta(K):=\max\{q_j:q_j\le K\}.
\end{equation*}
\begin{lemma}[General Denjoy--Koksma inequality]\label{lem2.3} Let $\theta\notin\Q$ and $0<t<1/2$. For any integers $K\ge 1$ and $u\ge 1$,
\[
\#\big\{u\le h<u+K: \|h\theta\|\le t\big\}\le 4Kt+\frac{4K}{Q_\theta(K)}.
\]
\end{lemma}

\subsection{Measure estimates}
The estimates below are largely based on results of \cite{HLW}, see Propositions 4.1, 4.2, 6.1, 6.2 and Lemma 11.1 there. To avoid the largely copy but try to make
 the current paper self-contained, we state the result in all cases but give a proof only for $n=m=1$ and $\theta\not\in \mathbb{Q}$.

\begin{proposition}\label{cor:measure2} For any $0<\eta<1$ and $T>1$, for all sufficiently large regular scales $M$, we have
\begin{itemize} \item[(1)] for $n=1$, \begin{equation}\label{count1}
\sum_{q\in\mathcal Q_M}\mathcal{L}_m\left(\mathcal{E}_{\eta}(q)\cap \bigcup_{s>T}D_{\eta}(s)\right)
\ll_m\left(\sum_{h>T}\psi(h)^m+ o_{\eta}(1)\right)\cdot\sum_{q\in\mathcal Q_M}\mathcal L_m(\mathcal E_{\eta}(q)).
\end{equation}
\item[(2)] for $n\ge 2$,
\begin{equation*}
\frac{1}{\mathcal{L}_{mn}(\mathcal{E}_{\eta}(\bq))}\sum_{|\bb|>T}\mathcal{L}_{mn}\left(\mathcal{E}_{\eta}(\bq)\cap\mathcal{D}_{\eta}(\bb)\right)
\ll_{m,n}\sum_{h>T}h^{n-1}\psi(h)^m+o_{\eta}(1).
\end{equation*}
\end{itemize}
\end{proposition}

\begin{proof}
(1) We show the case when $n=m=1$ and $\theta\notin \mathbb{Q}$ only. Fix $$0<\varepsilon<\frac16\left(\tau_{\psi}-\frac1m\right).$$
Let $M$ be a sufficiently large regular scales as in Lemma \ref{l1}. Suppose that
\begin{align*}
A_{\eta}(p,q)\cap D_{\eta}(r,s)\neq\varnothing.
\end{align*}
Then, there exists $k\in\mathbb{Z}$ such that
\begin{align}\label{ineq}
|k+(s-q)\theta|\le q\psi(s)+s\psi(q)\ \ \text{and}\ \ k\equiv sp\ (\text{mod }q),\  k\equiv -qr\ (\text{mod }s).
\end{align}

We distinguish the following two cases according to whether $q$ divides $s$.
\medskip

\noindent\textbf{Case 1: $q\nmid s$.} Since $q\in\mathcal Q_M$ is prime, we have $(q,s)=1$. Then there is a one-to-one correspondence between $k$ 
and the pair $(p\bmod q, \ r\bmod s)$.
We further divide the range of $s>T$ into the following four intervals, denoted by $I_1,I_2,I_3$ and $I_4$, respectively,
\begin{align*}
T< s<\frac{M}{\log M},\ \ \frac{M}{\log M}\le s<M,\ \  M\le s<M^{\tau_{\psi}-2\varepsilon},\ \ s\ge M^{\tau_{\psi}-2\varepsilon}.
\end{align*}

(a) For $s\in I_1\cup I_4$, it concludes
\begin{align*}
    &\sum_{q\in\mathcal Q_M}\sum_{\substack{s\in I_1\cup I_4}}\mathcal L_1(\mathcal E_{\eta}(q)\cap D_{\eta}(s))\\
    \le&\ \sum_{q\in\mathcal Q_M}\sum_{\substack{s\in I_1\cup I_4}}\#\{k\in\mathbb Z:\eqref{ineq}\text{ holds}\}\cdot\min\left\{\ell_{M}(\eta),\frac{\psi(s)}{s}\right\}\\
    \ll&\ \sum_{q\in\mathcal Q_M}\sum_{s\in I_1}(M\psi(s)+1)\cdot \ell_M(\eta)+\sum_{q\in\mathcal Q_M}\sum_{s\in I_4}(s\psi(q)+1)\cdot \frac{\psi(s)}{s}.
\end{align*}
The first sum on the right-hand side can be estimated directly as
$$\frac{M}{\log M}\cdot(1-\eta)\psi(M)\cdot\sum_{s>T}\psi(s)+\frac{M}{(\log M)^2}\cdot(1-\eta)\psi(M).$$
To estimate the second sum, we use Lemma \ref{l1}, which gives that for $q\in\mathcal Q_M$ and $s\in I_4$,
$$\psi(q)\ll M^{2\varepsilon}\psi(M),\ \ \psi(s)\le s^{-\tau_{\psi}+\varepsilon}\ \text{ and }\ \psi(s)\le \psi(M)\left(\frac{s}{M}\right)^{-\tau_{\psi}+\varepsilon}.$$
Substituting $\psi(q), \psi(s)$ by the above inequalities, the second sum is bounded by
$$\frac{M}{\log M}\cdot\left(M^{2\varepsilon}\psi(M)M^{(\tau_{\psi}-2\varepsilon)(1-\tau_{\psi}+\varepsilon)}+\psi(M)M^{(\tau_{\psi}-\varepsilon)(1-\tau_{\psi}+2\varepsilon)}\right).$$
Combined the above estimates for the two sums, we obtain the desired estimate \eqref{count1}.

(b) For $s\in I_2$, the inequality \eqref{ineq} implies that
$$|k+(s-q)\theta|\le 2M\left(\frac{M}{\log M}\right)^{-\tau_{\psi}+\varepsilon}<1/2,  \ \ {\text{and}}\ \
\Vert(s-q)\theta\Vert\le 2M\left(\frac{M}{\log M}\right)^{-\tau_{\psi}+\varepsilon}.$$
So there is at most one integer $k$ satisfying \eqref{ineq} for each pair $(q,s)$ and additionally $(q,s)$ also need to satisfy the second inequality above. We use the Denjoy-Koksma inequality (Lemma \ref{lem2.3}) to count the number of $s$ for a fixed $q$ which follows that
$$\#\left\{\frac{M}{\log M}<s<M: \|(s-q)\theta\|\le 2 M \bigg(\frac{M}{\log M}\bigg)^{-\tau_\psi+\varepsilon}\right\}\ll 2 M^2\bigg(\frac{M}{\log M}\bigg)^{-\tau_\psi+\varepsilon}+o(M).$$
Therefore, we have
\begin{align*}
    \sum_{q\in\mathcal Q_M}\sum_{s\in I_2}\mathcal L_1(\mathcal E_{\eta}(q)\cap D_{\eta}(s))&\ll\sum_{q\in\mathcal Q_M} \left(2M^2\bigg(\frac{M}{\log M}\bigg)^{-\tau_\psi+\varepsilon}+o(M)\right)\cdot\ell_M(\eta)\\
    &=o_{\eta}\left(\sum_{q\in\mathcal Q_M}\mathcal L_1(\mathcal E_{\eta}(q))\right).
\end{align*}

(c) For $s\in I_3$, similarly, \eqref{ineq} implies
$$|k+(s-q)\theta|\le 2^{1+\tau_{\psi}-\varepsilon}sM^{-\tau_{\psi}+\varepsilon}, \ \ \ \ \Vert(s-q)\theta\Vert\le 2^{1+\tau_{\psi}-\varepsilon}sM^{-\tau_{\psi}+\varepsilon}.$$
Since $s<M^{\tau_{\psi}-2\varepsilon}$, the right-hand side is $o(1)$. Hence, for sufficiently large $M$, there is at most one integer $k$ satisfying \eqref{ineq} for each pair $(q,s)$ and additionally $(q,s)$ also need to satisfy the second inequality above.

We partition the range of $s\in I_3$ into consecutive blocks of length at most $M$, each of form
$S\le s<S+M.$
For each block, we apply the Denjoy-Koksma inequality again, for $M$ sufficiently large, it follows that
\[
\#\Big\{S\le s<S+M:\|(s-q)\theta\|\le 2^{2+\tau_{\psi}-\varepsilon} SM^{-\tau_{\psi}+\varepsilon}\Big\}\ll SM^{1-\tau_{\psi}+\varepsilon}+o(M).
\]
Each corresponding intersection $A_\eta(p,q)\cap D_\eta(r,s)$ has length at most
\[
\ll \frac{\psi(s)}s\le\frac{\psi(S)}S.
\]
Moreover, Lemma \ref{l1} gives, for every integer $i\ge 1$,
$$\psi(iM)\le\psi(M)\left(\frac{1}{i}\right)^{\tau_{\psi}-\varepsilon}.$$
Combined this with the preceding estimates, we obtain
\begin{align*}
    \sum_{q\in\mathcal Q_M}\sum_{s\in I_3}\mathcal L_1(\mathcal E_{\eta}(q)\cap D_{\eta}(s))&\ll\sum_{q\in\mathcal Q_M}\sum_{1\le i\le M^{\tau_{\psi}-2\varepsilon-1}}\left(iM^{2-\tau_{\psi}+\varepsilon}+o(M)\right)\cdot\frac{\psi(iM)}{iM}\\
    &\ll\frac{M}{\log M}\cdot\psi(M)\left(M^{1-\tau_{\psi}+\varepsilon}+o(1)\right)\\
    &=o_{\eta}\left(\sum_{q\in\mathcal Q_M}\mathcal L_1(\mathcal E_{\eta}(q))\right).
\end{align*}

\noindent\textbf{Case 2: $q\mid s$.} At first, we claim that for every fixed $K$, for all sufficiently large $M$, and all $q\in\mathcal Q_M$, we have
\begin{equation*}
\mathcal E_\eta(q)\cap D_\eta(kq)=\emptyset, \ \ {\text{for all }}k\le K.
\end{equation*} This is clear true for $k=1$. Assume that for some $2\le k\le K$ and for some $x$ one has $$
\|qx-\theta\|<\psi(q), \ \|kqx-\theta\|<\eta \psi(kq).
$$ Then triangle inequality gives that$$
\|k\theta-\theta\|<k\psi(q)+\eta\psi(kq)<(1+K)\psi(q)
$$ which is impossible when $q$ is sufficiently large for irrational $\theta$.

On the other hand, for fixed $k$, the centers of $D_{\eta}(r,kq)$ are $1/(kq)$-separated, so each $A_\eta(p,q)$ can intersect at most
$$O\bigl(kq\ell_M(\eta)+1\bigr)$$
dangerous cubes $D_{\eta}(r,kq)$. Moreover, Lemma \ref{l1} gives that for integer $k>2$,
$$\psi(kq)\ll\psi(M)k^{-(\tau_{\psi}-\varepsilon)}.$$
Consequently, for each $q\in\mathcal Q_M$, for large $K$,
\begin{align*}
\sum_{k\ge 1}\mathcal L_1\left(\mathcal E_\eta(q)\cap D_\eta(kq)\right)
&=\sum_{k\ge K}\mathcal L_1\left(\mathcal E_\eta(q)\cap D_\eta(kq)\right)\\&\ll \sum_{k\ge K}q\cdot\bigl(kq\ell_M(\eta)+1\bigr)\cdot\frac{\psi(kq)}{kq}\\
&\ll\sum_{k\ge K}\left(k^{-(\tau_\psi-\varepsilon)}+k^{-(1+\tau_\psi-\varepsilon)}\right)\cdot\mathcal L_1(\mathcal E_{\eta}(q)).
\end{align*}
Since $\tau_{\psi}-\varepsilon>1$, the series on the right converges. Therefore, choose $K$ sufficiently large, we have
$$\sum_{q\in\mathcal Q_M}\sum_{k\ge 1}\mathcal L_1\left(\mathcal E_\eta(q)\cap D_\eta(kq)\right)\ll o_{\eta}\left(\sum_{q\in\mathcal Q_M}\mathcal L_1(\mathcal E_{\eta}(q))\right).$$

\end{proof}

\section{Auxiliary estimates}
\label{s:auxiliary}
Starting from Bluhm \cite{Bluhm98}, through a series works of Fraser and Hambrook \cite{FraH}, Fraser and Wheeler \cite{FraW}, Hambrook \cite{Ham, Ham1}, there has formed a standard method to estimate the Fourier dimension of a Diophantine set from below.

We first use smooth functions to approximate the characteristic function of the exact set and establish suitable Fourier decay for the resulting smooth approximating functions.
We then define a measure giving full mass to $\Exact_{m,n}(\psi,\bta)$ and estimate its Fourier coefficients to derive the lower bound
$$
\dim_{{\rm{F}}}\Exact_{m,n}(\psi,\bta)\ge\frac{2n}{1+\tau_\psi}.
$$

\subsection{Periodic bump function}
Choose a non-negative function $\varphi$ such that
\begin{itemize}
\item $\varphi$ is smooth, compactly supported on $(0,1)^m$, and satisfies
\begin{align}\label{bump_int}
\int_{\R^m}\varphi(y)\,\dif y=1;
\end{align}
\item for each integer $K\ge 1$, one has
\begin{equation}\label{bump_dec}|
\widehat\varphi(\xi)|\ll_K(1+|\xi|)^{-K},\quad \xi\in\R^m.	
\end{equation}
\end{itemize}

We define a periodic bump function as follows. Fix an integer \(M\ge 1\). For any \(0<\eta<1\) and any integer \(M/2<q\le M\),
define
\begin{equation*}
\phi_{M,q,\eta}(y)=\left(q\ell_M(\eta)\right)^{-m}\varphi\left(\frac{y-\eta\psi(q)\mathbf{1}}{q\ell_M(\eta)}\right),
\end{equation*}
where \(\mathbf{1}=(1,\ldots,1)\in\mathbb{R}^m\). It is clear that
\(\phi_{M,q,\eta}\) is non-negative and satisfies
\begin{equation*}
\supp \phi_{M,q,\eta}\subset\left[\eta\psi(q),\,\eta\psi(q)+q\ell_M(\eta)\right]^m.
\end{equation*}
We then define the periodic extension of \(\phi_{M,q,\eta}\) by \(\mathbb{Z}^m\)-periodization:
\[
\Phi_{M,q,\eta}(y)=\sum_{\bk\in\mathbb{Z}^m}\phi_{M,q,\eta}(y+\bk).
\]
Thus, \(\Phi_{M,q,\eta}\) is a well-defined function on the torus
\(\mathbb{T}^m\). A direct computation of its Fourier coefficients gives the following lemma.
\begin{lemma}\label{l31} For any $\bk\in\mathbb{Z}^m$, we have
$$\widehat{\Phi_{M,q,\eta}}(\bk)=\widehat{\phi_{M,q,\eta}}(\bk)=\widehat{\varphi}(q\ell_M(\eta)\cdot \bk)\cdot e^{-2\pi i\eta\psi(q)\bk}.$$
\end{lemma}

Recall that $\bta\in\mathbb{R}^m$. For any $\bq=(\bu,q)\in\mathcal{Q}_M$, define
\[
h_{M,\bq,\eta}(A)=\Phi_{M,q,\eta}(A\bq-\bta),\qquad A\in\Mat_{m,n}(\mathbb{R}).
\]
Since $\Phi_{M,q,\eta}$ is $\mathbb{Z}^m$-periodic and $\bq\in\mathbb{Z}^n$, the function $h_{M,\bq,\eta}$ is $\mathbb{Z}^{m\times n}$-periodic, and hence can be seen as a function on $\Mat_{m,n}(\mathbb{T})$. We state its basic properties in the following lemma.

\begin{lemma}\label{l32} For any $0<\eta<1$ and any $\bq=(\bu,q)\in\mathcal{Q}_M$, the following hold.
\begin{enumerate}[\upshape (1)]
\item For every $A\in\Mat_{m,n}(\mathbb{T})$, we have
\[
0\le h_{M,\bq,\eta}(A)\ll \frac{\mathbf{1}_{\mathcal{E}_{\eta}(\bq)}(A)}{\mathcal L_{mn}(\mathcal E_\eta(\bq))}.
\]
In particular, $\supp h_{M,\bq,\eta}\subset\mathcal{E}_{\eta}(\bq).$

\item The function $h_{M,\bq,\eta}$ has integral $1$, that is,
\[
\int_{\Mat_{m,n}(\mathbb{T})} h_{M,\bq,\eta}(A)\,\dif A=1.
\]

\item For every $\xi\in\mathbb{Z}^{m\times n}$,
\[
\widehat{h_{M,\bq,\eta}}(\xi)=
\begin{cases}
\widehat{\Phi_{M,q,\eta}}(\bk) e^{-2\pi i\bk\cdot\bta},&\text{if } \xi=\bk\bq^\top\text{ for some }\bk\in\mathbb{Z}^m, \\[4pt]
0,&\text{otherwise}.
\end{cases}
\]
Consequently, for every integer $K\ge1$,
\[
\left|\widehat{h_{M,\bq,\eta}}(\xi)\right|\ll_K\left(1+(1-\eta)\frac{\psi(M)}{M}|\xi|\right)^{-K}
\]
whenever $\xi=\bk\bq^\top$ for some $\bk\in\mathbb{Z}^m$.
\end{enumerate}
\end{lemma}
\begin{proof} \textup{(1)} Since $\phi_{M,q,\eta}$ is non-negative and
\[
\supp \phi_{M,q,\eta} \subset\left[\eta\psi(q),\,\eta\psi(q)+q\ell_M(\eta)\right]^m,
\]
we have
\[
h_{M,\bq,\eta}(A)\ge0 \quad \text{and} \quad \supp h_{M,\bq,\eta}\subset\mathcal{E}_{\eta}(\bq).
\]
Moreover, since $M/2<q\le M$, by the definition of $\ell_M(\eta)$ we have
\begin{equation}\label{supp_h}
|\phi_{M,q,\eta}|\ll (q\ell_M(\eta))^{-m}\stackrel{\eqref{measure}}{=}\frac{1}{\mathcal L_{mn}(\mathcal E_\eta(\bq))}.
\end{equation}
This together with \eqref{supp_h} gives
\[
0\le h_{M,\bq,\eta}(A)\ll\frac{\mathbf{1}_{\mathcal{E}_{\eta}(\bq)}(A)}{\mathcal L_{mn}(\mathcal E_\eta(\bq))}.
\]

\textup{(2)} Since $\bq\neq\mathbf{0}$, the continuous group homomorphism
\[
\pi_{\bq}:\Mat_{m,n}(\mathbb{T})\to\mathbb{T}^m, \quad \pi_{\bq}(A)=A\bq,
\]
is surjective. Hence $\pi_{\bq}$ pushes forward the normalized Haar measure on $\Mat_{m,n}(\mathbb{T})$ to normalized Haar measure on $\mathbb{T}^m$. By translation invariance of Haar measure,
\[
\begin{aligned}
\int_{\Mat_{m,n}(\mathbb{T})}h_{M,\bq,\eta}(A)\,\dif A &=\int_{\Mat_{m,n}(\mathbb{T})}\Phi_{M,q,\eta}(A\bq-\bta)\,\dif A\\
&=\int_{\mathbb{T}^m}\Phi_{M,q,\eta}(y-\bta)\,\dif y\\
&=\int_{\mathbb{T}^m}\Phi_{M,q,\eta}(y)\,\dif y=1,
\end{aligned}
\]
where the last equality follows from \eqref{bump_int}.

\textup{(3)} Expanding $\Phi_{M,q,\eta}$ into its Fourier series gives
\begin{align*}
h_{M,\bq,\eta}(A)&=\Phi_{M,q,\eta}(A\bq-\bta)=\sum_{\bk\in\mathbb{Z}^m}\widehat{\Phi_{M,q,\eta}}(\bk)e^{2\pi i\bk\cdot(A\bq-\bta)}\\
&=\sum_{\bk\in\mathbb{Z}^m}\widehat{\Phi_{M,q,\eta}}(\bk)e^{-2\pi i\bk\cdot\bta}e^{2\pi i\langle \bk\bq^\top,A\rangle}.
\end{align*}
Therefore,
\[
\widehat{h_{M,\bq,\eta}}(\xi)=\begin{cases}
\widehat{\Phi_{M,q,\eta}}(\bk) e^{-2\pi i\bk\cdot\bta}, & \text{if } \xi=\bk\bq^\top \text{ for some }\bk\in\mathbb{Z}^m, \\[4pt]
0, & \text{otherwise}.
\end{cases}
\]
Finally, by Lemma~\ref{l31} and \eqref{bump_dec}, if $\xi=\bk\bq^\top$, then
\begin{align*}
\big|\widehat{h_{M,\bq,\eta}}(\xi)\big|
&=\big|\widehat{\Phi_{M,q,\eta}}(\bk)\big|=\big|\widehat{\varphi}\bigl(q\ell_M(\eta)\bk\bigr)\big|\ll_K\left(1+q\ell_M(\eta)|\bk|\right)^{-K}.
\end{align*}
Since $|\xi|=q|\bk|$ and $q\ell_M(\eta)=q\ell_M(\eta)$, we have
\[
q\ell_M(\eta)|\bk|=\ell_M(\eta)|\xi|=(1-\eta)\frac{\psi(M)}{M}|\xi|.
\]
Thus,
\[
\big|\widehat{h_{M,\bq,\eta}}(\xi)\big|\ll_K\left(1+(1-\eta)\frac{\psi(M)}{M}|\xi|\right)^{-K},
\]
which completes the proof.
\end{proof}

\subsection{Single-scale function}
For any $M\ge1$ and $A\in\Mat_{m,n}(\mathbb{T})$, define
\[
f_{M,\eta}(A)=\frac{1}{\#\mathcal{Q}_M}\sum_{\bq\in\mathcal{Q}_M}h_{M,\bq,\eta}(A).
\]
The following properties follow immediately from the definition of $f_{M,\eta}$ and Lemma~\ref{l32}.

\begin{lemma}\label{l33} 
The function $f_{M,\eta}$ satisfies:
\begin{enumerate}[\upshape (1)]
\item $f_{M,\eta}$ is non-negative and
\[\int_{\Mat_{m,n}(\mathbb{T})}f_{M,\eta}(A)\,\dif A=1.
\]

\item If $f_{M,\eta}(A)\neq0$, then there exists $\bq\in\mathcal{Q}_M$ such that
\[
\eta\psi(|\bq|)\le\Vert A\bq-\bta\Vert<\psi(|\bq|).
\]
In particular, for any sequences $M_j\to\infty$ and $\eta_j\in(0,1)$,
\[
\bigcap_{j=1}^{\infty} \supp(f_{M_j,\eta_j})\subset W_{m,n}(\psi,\bta).
\]
\end{enumerate}
\end{lemma}
\begin{remark}\label{r2} The measure $\mu$ to be constructed later will be supported on $\cap_{j=1}^{\infty}\supp(f_{M_j,\eta_j})$ but gives full measure to $\Exact_{m,n}(\psi,\bta)$. So we are in a very similar setting as estimating the Fourier dimension of $W_{m,n}(\psi,\bta)$.
\end{remark}
For any $0<\ve<\tau_\psi-n/m$, write
\begin{equation}\label{eq:beta-def}
\beta=\beta(\tau_\psi,\ve)=1+\tau_{\psi}+2\varepsilon.
\end{equation}

\begin{lemma}\label{l34} For any $0<\eta<1$ and sufficiently large $M$, the following hold.
\begin{enumerate}[\upshape (1)]
\item For any $0<|\xi|<M/2,$ we have $\widehat {f_{M,\eta}}(\xi)=0;$
		
\item For any $0<|\xi|\le M^{\beta}$, we have $|\widehat {f_{M,\eta}}(\xi)|\ll_{\tau_\psi,\ve}M^{-n+\ve};$
			
\item For every integer $K\ge1$,
\begin{equation*}
|\widehat{f_{M,\eta}}(\xi)|\ll_K\left(1+(1-\eta)\frac{\psi(M)}{M}|\xi|\right)^{-K},\ \xi\in\Z^{m\times n}.
\end{equation*}
\end{enumerate}
\end{lemma}
\begin{proof} By the definition of $f_{M,\eta}$,
\[
\widehat{f_{M,\eta}}(\xi)=\frac{1}{\#\mathcal{Q}_M}\sum_{\bq\in\mathcal{Q}_M}\widehat{h_{M,\bq,\eta}}(\xi).
\]

\textup{(1)} By Lemma~\ref{l32},
$\widehat{h_{M,\bq,\eta}}(\xi)\neq0$ only if $\xi=\bk\bq^\top$ for some $\bk\in\mathbb{Z}^m$. If $\xi\neq0$, then $\bk\neq0$, and hence we have
$$|\xi|=|\bk|\,|\bq|\ge|\bq|>\frac{M}{2}.$$
Thus, if $0<|\xi|<M/2$, every term in the above sum vanishes, and $\widehat{f_{M,\eta}}(\xi)=0.$

\textup{(2)} We distinguish the cases $n\ge2$ and $n=1$.
\begin{itemize}
\item When $n\ge2$. We claim that, for each $\xi\neq0$, there is at most one $\bq\in\mathcal{Q}_M$ such that $\widehat{h_{M,\bq,\eta}}(\xi)\neq0$. Indeed, suppose that
\[
\xi=\bt\bq^\top=\bt'\bq'^\top\neq0
\]
for some $\bq,\bq'\in\mathcal{Q}_M$. Since $\xi\neq0$, it has a non-zero row, say the $i$-th row. Thus
\[
t_i\bq^\top =t_i'\bq'^\top\neq0,
\]
which implies that $\bq$ and $\bq'$ are parallel. Since both are primitive integer vectors, we must have $\bq'=\pm\bq.$ As their last coordinates are positive, it follows that $\bq'=\bq$. Hence,
\[
\left|\widehat{f_{M,\eta}}(\xi)\right|\le\frac{1}{\#\mathcal{Q}_M}\ll\frac{\log M}{M^n}\ll_{\varepsilon} M^{-n+\varepsilon}.
\]

\item When $n=1$. Let $0<|\xi|\le M^\beta$, and choose a non-zero coordinate $\xi_i$. If $\ell$ distinct primes $q\in(M/2,M]$ satisfy
    $\widehat{h_{M,q,\eta}}(\xi)\neq0$, then each such $q$ divides $\xi_i$. Therefore,
\[\left(\frac{M}{2}\right)^\ell\le |\xi_i| \le |\xi| \le M^\beta,
\]
and hence $\ell=O(1)$.
It follows that
\[
\left|\widehat{f_{M,\eta}}(\xi)\right|\le \ell \frac{\log M}{M} \ll M^{-1+\varepsilon}.
\]
\end{itemize}
This completes the proof of \textup{(2)}.

(3) Lemma~\ref{l32} gives, for every $\bq\in\mathcal{Q}_M$,
\[
\left|\widehat{h_{M,\bq,\eta}}(\xi)\right|\ll_K\left(1+{(1-\eta)}\frac{\psi(M)}{M}|\xi|\right)^{-K}.
\]
Averaging over $\bq\in\mathcal{Q}_M$ yields
\[
\left|\widehat{f_{M,\eta}}(\xi)\right|\ll_K\left(1+{(1-\eta)}\frac{\psi(M)}{M}|\xi|\right)^{-K},
\]
which proves \textup{(3)}.
\end{proof}

\subsection{Key estimates}
The following proposition shows that $W_{m,n}(\eta\psi,\bta)$ carries a small mass with respect to the measure \(f_{M,\eta}(A)\,\dif A\).
This is the most essential fact ensuring that the resulting measure we will construct is supported on $W_{m,n}(\psi,\bta)$ but gives
full mass to $\Exact_{m,n}(\psi,\bta)$.

\begin{proposition}\label{p31}
Let $0<\eta<1$ and $T>1$. Then, as $M\to\infty$ through regular scales, we have
\begin{equation}\label{bridge}
\int f_{M,\eta}(A)\mathbf{1}_{\bigcup_{|\bb|>T}\mathcal{D}_{\eta}(\bb)}(A)\,\dif A\le c_{m,n}\sum_{h>T} h^{n-1}\psi(h)^m+o_\eta(1),
\end{equation}
where $c_{m,n}$ is an absolute constant.
\end{proposition}
\begin{proof}
By the definition of $f_{M,\eta}$ and Lemma~\ref{l32}, we have
\begin{align*}
&\int f_{M,\eta}(A) \mathbf{1}_{\bigcup_{|\bb|>T}\mathcal{D}_{\eta}(\bb)}(A)\,\dif A\\
=&\ \frac{1}{\#\mathcal{Q}_M} \sum_{\bq\in\mathcal{Q}_M}\int
h_{M,\bq,\eta}(A)\mathbf{1}_{\bigcup_{|\bb|>T}\mathcal{D}_{\eta}(\bb)}(A)\,\dif A\\
\ll&\ \frac{1}{\#\mathcal{Q}_M}\sum_{\bq\in\mathcal{Q}_M}\frac{1}{\mathcal{L}_{mn}(\mathcal{E}_{\eta}(\bq))}
\int\mathbf{1}_{\mathcal{E}_{\eta}(\bq)}(A)\mathbf{1}_{\bigcup_{|\bb|>T}\mathcal{D}_{\eta}(\bb)}(A)\,\dif A\\
\le &\ \frac{1}{\#\mathcal{Q}_M}\sum_{\bq\in\mathcal{Q}_M}\frac{1}{\mathcal{L}_{mn}(\mathcal{E}_{\eta}(\bq))}\cdot\mathcal{L}_{mn}\left(\mathcal{E}_{\eta}(\bq)\cap\bigcup_{|\bb|>T}\mathcal{D}_{\eta}(\bb)\right).
\end{align*}
Applying Propostion~\ref{cor:measure2}, we obtain
\[
\int f_{M,\eta}(A)\mathbf{1}_{\bigcup_{|\bb|>T}\mathcal{D}_{\eta}(\bb)}(A)\,\dif A\ll_{m,n}\sum_{h>T}h^{n-1}\psi(h)^m+o_{\eta}(1),
\]
which proves \eqref{bridge}.
\end{proof}

\section{Construction of the measure}\label{s:con}
In this section, we construct a finite positive measure $\mu$ with polynomial Fourier decay and show that $\Exact_{m,n}(\psi,\bta)$ has full $\mu$-measure.

\subsection{Choice of parameters}
Fix some parameters first. Fix $0<s<{2n}/{(1+\tau_\psi)}$ and then choose $\alpha$ such that
\[
\frac{s}{2}<\alpha<\frac{n}{1+\tau_\psi}.
\]
Let $\ve>0$ be sufficiently small so that
\begin{equation*}
m(\tau_\psi-\ve)>n \quad\text{and}\quad \alpha(1+\tau_\psi+2\ve)+2\ve<n.
\end{equation*}
For ease of notation, recall from \eqref{eq:beta-def} that
\begin{equation}\label{beta}
\beta=1+\tau_\psi+2\ve\ \Longrightarrow\ \alpha\beta+2\ve<n.
\end{equation}
Finally, let $\eta_j=1-2^{-j-2}$ for each $j\ge1.$

\subsection{A stability lemma}
The following stability lemma allows us to choose the scales recursively without losing the Fourier decay obtained at the previous stages.

\begin{lemma}\label{l41}
Let $g$ be a function on $\Mat_{m,n}(\T)$ with $\widehat g(0)=1$, and suppose that
\[
\sup_{\xi\in\Z^{m\times n}\setminus\{0\}}(1+|\xi|)^\alpha |\widehat g(\xi)|\le \delta
\]
for some $\delta>0$. Then the following statements hold.

\begin{enumerate}[\upshape (1)]
\item For every smooth function $f$ on $\Mat_{m,n}(\T)$,
\[
\sup_{\xi\in\Z^{m\times n}}(1+|\xi|)^\alpha\big|\widehat{fg}(\xi)-\widehat f(\xi)\big|\le\delta\sum_{\zeta\in\Z^{m\times n}}(1+|\zeta|)^\alpha|\widehat f(\zeta)|.
\]

\item For every smooth function $f$ and every smooth function
$P$ on $\Mat_{m,n}(\T)$,
\[
\left|\int Pfg\,\dif A-\int Pf\,\dif A\right|\le\delta\left(\sum_{\xi\in\Z^{m\times n}}(1+|\xi|)^{-\alpha}|\widehat P(\xi)|\right)
\left(\sum_{\zeta\in\Z^{m\times n}}(1+|\zeta|)^\alpha|\widehat f(\zeta)|\right).
\]
\end{enumerate}
\end{lemma}
\begin{proof}
(1) Since $\widehat g(0)=1$, the convolution formula gives
\[
\widehat{fg}(\xi)-\widehat f(\xi)=\sum_{\zeta\ne0}\widehat f(\xi-\zeta)\widehat g(\zeta).
\]
Hence, by using the following inequality
$$1+|\xi|\le(1+|\xi-\zeta|)(1+|\zeta|),$$
it follows that
\begin{align*}
(1+|\xi|)^\alpha\big|\widehat{fg}(\xi)-\widehat f(\xi)\big|&\le \sum_{\zeta\ne0}\Big[(1+|\xi-\zeta|)^\alpha|\widehat f(\xi-\zeta)|\Big] \cdot \Big[(1+|\zeta|)^\alpha |\widehat g(\zeta)|\Big]\\ &\le\delta\sum_{\zeta\ne0}(1+|\xi-\zeta|)^\alpha|\widehat f(\xi-\zeta)|
\le\delta\sum_{\nu\in\Z^{m\times n}}(1+|\nu|)^\alpha|\widehat f(\nu)|.
\end{align*}

(2) Representing the smooth functions by their Fourier series, one has
\begin{align*}
\int Pfg\,\dif A-\int Pf\,\dif A&=\sum_{\xi\in\Z^{m\times n}}\widehat P(-\xi)\left(\widehat{fg}(\xi)-\widehat f(\xi)\right).
\end{align*}
Therefore, by~{\rm (1)}, we have
\begin{align*}
\left|\int Pfg\,\dif A-\int Pf\,\dif A\right|&\le\sum_{\xi}|\widehat P(-\xi)|\big|\widehat{fg}(\xi)-\widehat f(\xi)\big|\\
&\le\delta\left(\sum_{\xi}(1+|\xi|)^{-\alpha}|\widehat P(\xi)|\right)\left(\sum_{\zeta}(1+|\zeta|)^\alpha|\widehat f(\zeta)|
\right).\qedhere
\end{align*}
\end{proof}
We will apply the stability lemma to functions \(f_{M,\eta}\), so we show that they satisfy the required bound in the above lemma.

\begin{lemma}\label{l42}
There exists $\kappa>0$ such that, for every fixed $j\ge1$ and all sufficiently large regular scales $M$,
\[
\sup_{\xi\in\Z^{m\times n}\setminus\{0\}}(1+|\xi|)^\alpha|\widehat f_{M,\eta_j}(\xi)|\le M^{-\kappa}.
\]
\end{lemma}

\begin{proof}
Fix $j\ge1$. Recall \eqref{beta} that $$
\beta=1+\tau_\psi+2\ve,\ \ {\text{and}} \ \  \alpha\beta+2\ve<n.
$$ By Lemma~\ref{l34}\textup{(1)--(2)}, for $0<|\xi|\le M^\beta$,
\[
(1+|\xi|)^\alpha|\widehat f_{M,\eta_j}(\xi)|\ll M^{\alpha\beta-n+\ve}\ll M^{-\ve}.
\]
For $|\xi|>M^\beta$,  by Lemma~\ref{l34}\textup{(3)} and $\psi(M)\ge M^{-(\tau_\psi+\ve)}$, we have
\begin{align*}
(1+|\xi|)^\alpha|\widehat f_{M,\eta_j}(\xi)|&\ll_{K,\alpha,\eta_j}(1+|\xi|^{\alpha})\cdot \left(1+M^{-(1+\tau_{\psi}+\epsilon)}|\xi|\right)^{-K}\\
&\asymp |\xi|^{\alpha}\cdot \left(M^{-(1+\tau_{\psi}+\epsilon)}|\xi|\right)^{-K}
\le M^{(\alpha-K)\beta}\cdot M^{K(1+\tau_\psi+\ve)}\\
&= M^{\alpha\beta-K\ve}.
\end{align*} 
Choosing $K$ sufficiently large such that $K\ve>\alpha\beta+\ve$, gives
\[
(1+|\xi|)^\alpha|\widehat f_{M,\eta_j}(\xi)|\ll_{K,\alpha,\eta_j}M^{-\ve}.
\]
Consequently, for every fixed $j$ and all sufficiently large regular
scales $M$,
\[
\sup_{\xi\ne0}(1+|\xi|)^\alpha|\widehat f_{M,\eta_j}(\xi)|\le M^{-\ve/2}.\qedhere
\]
\end{proof}

\subsection{Measure Construction and its Fourier decay}

We now construct the measure inductively. We shall choose a strictly increasing sequence of regular scales $M_1<M_2<\cdots.$ Once $M_j$ is chosen, we write $f_j=f_{M_j,\eta_j}$ and define
\[
\dif\mu_j=\prod_{i=1}^j f_i\,\dif A.
\]
For convenience, let $\mu_0=\mathcal L_{mn}.$ Along with the scales $M_j$, we choose an increasing sequence $\{T_j\}_{j\ge1}$ and a sequence of smooth functions $\{P_j\}_{j\ge1}$. These sequences are chosen to satisfy the following properties:
\begin{enumerate}[\upshape (a)]
\item For every $j\ge1$,
\[
\prod_{i=1}^{j-1}\|f_i\|_\infty\int f_j(A)\mathbf1_{\bigcup_{|\bb|>T_j}\mathcal D_{\eta_j}(\bb)}(A)\,\dif A <2^{-j-2},
 \]
where the empty product is understood to be $1$ when $j=1$.

\item For every \(j\ge1\),
\begin{equation}\label{eq1}
\sup_{\xi\in\Z^{m\times n}}(1+|\xi|)^\alpha\left|\widehat\mu_j(\xi)-\widehat\mu_{j-1}(\xi)\right|<2^{-j-2}.
\end{equation}

\item For every $j\ge1$, let
\begin{equation}\label{eq3}\mathcal B_j=\bigcup_{\substack{\bb\in\Z^n\setminus\{\bzero\}\\T_j<|\bb|\le T_{j+1}}}\mathcal D_{\eta_j}(\bb).
\end{equation}
The smooth function $P_j$ approximates the characteristic function of $\mathbf1_{\mathcal B_j}$ satisfying
\[
P_j\ge\mathbf1_{\mathcal B_j}, \quad \int P_j\,\dif\mu_j\le\mu_j(\mathcal B_j)+2^{-j-2}.
\] and for every \(j\ge2\) and every \(1\le k<j\),
\begin{equation}\label{eq2}
\left|\int P_k\dif\mu_j-\int P_k\dif\mu_{j-1}\right|<2^{-j-2}.
\end{equation}
\end{enumerate}
\begin{remark}
The above quantities are chosen recursively in the following order:
\[
T_1
\longrightarrow M_1 \longrightarrow T_2 \longrightarrow P_1 \longrightarrow M_2 \longrightarrow T_3 \longrightarrow P_2 \longrightarrow M_3
\longrightarrow\cdots.
\]
More precisely, once $T_j$ has been chosen, we choose $M_j$ sufficiently large so that \textup{(a)} and \textup{(b)} hold. Then choose $T_{j+1}$ sufficiently large, both for the next stage and to determine the set $\mathcal B_j$. At last choose $P_j$ so that \textup{(c)} holds.
\end{remark}

We now verify that the above construction can be carried out. Recall that
$$\sum_{h=1}^{\infty}h^{n-1}\psi(h)^m<\infty,\ \ {\text{since}}\ \tau_{\psi}>n/m$$
so the right series in Proposition~\ref{p31} can be arbitrarily small if $T$ is large and $M$ is large accordingly.

Let us start with $j=1$. Choose $T_1>1$ sufficiently large and then choose a sufficiently large regular scale $M_1$ so that, by Proposition~\ref{p31},
\[
\int f_1(A)\mathbf1_{\bigcup_{|\bb|>T_1}\mathcal D_{\eta_1}(\bb)}(A)\,\dif A<2^{-3}.
\]
Thus \textup{(a)} holds for $j=1$.

We can increase $M_1$, if necessary, so that \textup{(b)} also holds. Indeed, since
\[
\dif\mu_1=f_1\,\dif A, \quad \mu_0=\mathcal L_{mn},
\]
we have $\widehat\mu_1(0)=\widehat\mu_0(0)=1,$
and for $\xi\neq0$,
\[
\widehat\mu_1(\xi)-\widehat\mu_0(\xi)=\widehat\mu_1(\xi)=\widehat f_1(\xi).
\]
By Lemma \ref{l42}, we can further ask the regular scale $M_1$ so large that
\[
\sup_{\xi\in\Z^{m\times n}}(1+|\xi|)^\alpha\left|\widehat\mu_1(\xi)-\widehat\mu_0(\xi)\right|=\sup_{\xi\in\Z^{m\times n}\setminus\{0\}}(1+|\xi|)^\alpha|\widehat f_{1}(\xi)|\le M^{-\kappa}<2^{-3}.
\]

Having fixed $M_1$, choose $T_2>T_1$ sufficiently large so that
\[c_{m,n}\|f_1\|_{\infty}\sum_{h>T_2} h^{n-1}\psi(h)^m<2^{-2-3}.\]
Let $\mathcal B_1$ be defined as in \eqref{eq3}. Since
\[
\mathcal B_1\subseteq \bigcup_{|\bb|>T_1}\mathcal D_{\eta_1}(\bb),
\]
property \textup{(a)} gives
$$\mu_1(\mathcal B_1)=\int f_1(A)\mathbf1_{\mathcal B_j}\dif A\le \int f_1(A)\mathbf1_{\bigcup_{|\bb|>T_1}\mathcal D_{\eta_1}(\bb)}(A)\,\dif A <2^{-3}.$$

For any nonzero $\bb\in\Z^n$, the map $A\mapsto A\bb-\bta$ pushes Haar measure on $\Mat_{m,n}(\T)$ forward to Haar measure on $\T^m$. The boundary of a sup-norm ball in $\T^m$ has Haar measure zero. Hence $\cL_{mn}(\partial\cd_{\eta_1}(\bb))=0$; because $\mathcal B_1$ is a finite union, this implies $\cL_{mn}(\partial\mathcal B_1)=0$. Since $\dif\mu_1=f_1\dif A$,
\[
\mu_1(\overline{\mathcal B_1})=\mu_1(\mathcal B_1).
\]
By outer regularity we may choose an open neighborhood $V_1\supset\overline{\mathcal B_1}$ with $\mu_1(V_1)<\mu_1(\mathcal B_1)+2^{-3}$. By the existence of smooth bump functions (see \cite[Proposition~2.25]{Lee}), there exists $P_1\in C^\infty(\Mat_{m,n}(\T))$ such that
\[
0\le P_1\le1,\qquad P_1=1\ \text{on }\overline{\mathcal B_1},\qquad \supp P_1\subset V_1.
\]
Consequently,
\[
P_1\ge\mathbf1_{\mathcal B_1},\qquad \int P_1\,\dif\mu_1\le\mu_1(V_1)<\mu_1(\mathcal B_1)+2^{-3}.
\]
This establishes \textup{(c)} at stage $1$.

We next choose $M_2$. Since $T_2$ has been given, choose $M_2$ large enough so that \textup{(a)} holds for $j=2$. Since $f_1$ and $P_1$ are smooth functions, one has
$$
\left(\sum_{\xi}(1+|\xi|)^{-\alpha}|\widehat P_1(\xi)|\right)<\infty, \ \ \ \left(\sum_{\zeta}(1+|\zeta|)^\alpha|\widehat f_1(\zeta)|\right)<\infty.
$$
Then ask $M_2$ so large that \begin{align*}
\sup_{\xi\in\Z^{m\times n}}(1+|\xi|)^\alpha\left|\widehat\mu_2(\xi)-\widehat\mu_1(\xi)\right|&=\sup_{\xi\in\Z^{m\times n}}(1+|\xi|)^\alpha\big|\widehat{f_1f_2}(\xi)-\widehat f_1(\xi)\big|\\
&\le M_2^{-\kappa}
 \left(\sum_{\zeta}(1+|\zeta|)^\alpha|\widehat f_1(\zeta)|\right)<2^{-4}.
\end{align*}
Finally, applying Lemma~\ref{l41}\textup{(2)} with $f=f_1$, $g=f_2$, $P=P_1,$ and ask $M_2$ sufficiently large,
we obtain
\begin{align*}
&\left|\int P_1\,\dif\mu_2-\int P_1\,\dif\mu_1\right|\le \ M_2^{-\kappa}\left(\sum_{\xi}(1+|\xi|)^{-\alpha}|\widehat P_1(\xi)|\right)\left(\sum_{\zeta}(1+|\zeta|)^\alpha|\widehat f_1(\zeta)|\right)<2^{-4}.
\end{align*}
Thus \textup{(a)} and \textup{(b)} hold for $j=2$.

We have now constructed
\[
T_1\longrightarrow M_1\longrightarrow T_2\longrightarrow P_1\longrightarrow M_2.
\]
We proceed inductively. Suppose that, for some $j\ge2$, we have constructed
\[
T_1\longrightarrow M_1\longrightarrow T_2\longrightarrow P_1\longrightarrow M_2\longrightarrow\cdots\longrightarrow T_j\longrightarrow P_{j-1}\longrightarrow M_j
\]
so that the required properties (a)--(c) hold. We now construct
\[
T_{j+1}\longrightarrow P_j\longrightarrow M_{j+1}.
\]

Since $f_1,\ldots,f_j$ have already been fixed and $\tau_\psi>{n}/{m}$, we first choose $T_{j+1}>T_j$ sufficiently large so that
\[
c_{m,n}\prod_{i=1}^{j}\|f_i\|_\infty\sum_{h>T_{j+1}}h^{n-1}\psi(h)^m<2^{-j-4}.
\]

For such $T_{j+1}$, let $\mathcal B_j$ be defined as in \eqref{eq3}. Then, \begin{align*}
  \mu_j(\mathcal B_j)=\int \prod_{i=1}^{j}f_i(A)\mathbf1_{\mathcal B_j}\dif A\le \prod_{i=1}^{j-1}\|f_i\|_\infty\int f_j(A)\mathbf1_{\bigcup_{|\bb|>T_j}\mathcal D_{\eta_1}(\bb)}(A)\,\dif A.
\end{align*} So Proposition~\ref{p31} is used to conclude that
\[
\mu_j(\mathcal B_j)<2^{-j-2}.
\]
Since $\mathcal B_j$ is a finite union, the same smooth approximation argument as in the case $j=1$ allows us to choose $P_j$ such that
\[
P_j\ge\mathbf1_{\mathcal B_j}, \quad \int P_j\,\dif\mu_j \le \mu_j(\mathcal B_j)+2^{-j-2}.
\]
Thus \textup{(c)} holds for $P_j$.

Finally, after \(T_{j+1}\) and \(P_j\) have been fixed, we choose \(M_{j+1}>M_j\) to be a sufficiently large regular scale. First, Proposition~\ref{p31}, together with the above choice of \(T_{j+1}\), ensures \textup{(a)} at stage \(j+1\) once \(M_{j+1}\) is sufficiently large. We then increase \(M_{j+1}\), if necessary, so that Lemma~\ref{l42} and Lemma~\ref{l41}\textup{(1)} give
\begin{align*}
\sup_{\xi\in\Z^{m\times n}}(1+|\xi|)^\alpha\left|\widehat\mu_{j+1}(\xi)-\widehat\mu_j(\xi)\right|&=\sup_{\xi\in\Z^{m\times n}}(1+|\xi|)^\alpha\big|\widehat{f_1f_2\cdots f_{j+1}}(\xi)-\widehat {f_1\cdots f_{j}}(\xi)\big|\\
&\le M_{j+1}^{-\kappa}
 \left(\sum_{\zeta}(1+|\zeta|)^\alpha|\widehat {f_1\cdots f_{j}}(\zeta)|\right)<2^{-(j+1)-2}.
\end{align*}
which is \eqref{eq1} at stage \(j+1\). Moreover, since \(P_1,\ldots,P_j\) have already been fixed, Lemma~\ref{l41}\textup{(2)} allows us, by increasing \(M_{j+1}\) further if necessary, to ensure that
$$
\left|\int P_k\,\dif\mu_{j+1}-\int P_k\,\dif\mu_j\right|<2^{-(j+1)-2},\qquad 1\le k\le j.
$$
This is precisely \eqref{eq2} at stage \(j+1\). Since there are only finitely many smooth functions \(P_1,\ldots,P_j\), all these estimates can be satisfied simultaneously. Thus \textup{(b)} holds at stage \(j+1\), completing the inductive step.

Let $\mu$ be the weak limit of $\{\mu_j\}_{j\ge1}$. We will show that $\mu$ has the desired Fourier decay and gives full mass to $\Exact_{m,n}(\psi,\bta)$. 
\begin{remark}We give some explanation about the second claim which looks a little strange. In fact, $\mu$ is supported on
\[
\bigcap_{j=1}^{\infty} \supp(f_{M_j,\eta_j})\subset \bigcap_{j=1}^{\infty} \bigcup_{\bq\in\cq_{M_j}}\ce_{\eta_j}(\bq).
\]
The reason that $\mu$ gives full mass to $\Exact_{m,n}(\psi,\bta)$ is because the dangerous part
$\{\mathcal D_{\eta_j}(\bb): |\bb|>T_j\}$
takes only a small part in $\{\ce_{\eta_j}(\bq):\bq\in\cq_{M_j}\}$ as Proposition~\ref{p31} indicates.\end{remark}

\begin{proposition}
 $\mu$ is a non-zero finite measure and satisfies
\begin{enumerate}[\upshape (1)]
 \item $|\widehat\mu(\xi)|\ll (1+|\xi|)^{-\alpha}$, for any $\xi\in\Z^{m\times n}$.

\item $\mu\bigl(\Exact_{m,n}(\psi,\bta)\bigr) =\mu\bigl(\Mat_{m,n}(\T)\bigr).$
\end{enumerate}
Consequently,
\[
\fdim\Exact_{m,n}(\psi,\bta)\ge 2\alpha.
\]
Since $\alpha$ can be chosen arbitrarily close to $n/(1+\tau_\psi)$, we conclude
\[
\fdim\Exact_{m,n}(\psi,\bta)\ge\frac{2n}{1+\tau_\psi}.
\]
\end{proposition}

\begin{proof}
\textup{(1)}. By the construction of the measures $\{\mu_j\}$,
\[
\sup_{\xi\in\Z^{m\times n}}(1+|\xi|)^\alpha\left|\widehat\mu_j(\xi)-\widehat\mu_{j-1}(\xi)\right|<2^{-j-2}.
\]
Taking $\xi=0$ and summing over $j$, we obtain
\[
\left|\widehat\mu_j(0)-1\right|\le\sum_{i=1}^j2^{-i-2}<\frac14.
\]
Thus the total mass of $\mu_j$ are uniformly bounded and bounded away from zero. In particular, the weak limit $\mu$ is a non-zero finite measure.

For $\xi\neq0$, since $\widehat\mu_0(\xi)=0,$
we have
\begin{align*}
|\widehat\mu(\xi)|\le\sum_{j=1}^\infty\left|\widehat\mu_j(\xi)-\widehat\mu_{j-1}(\xi)\right|
\le(1+|\xi|)^{-\alpha}\sum_{j=1}^\infty2^{-j-2}\ll(1+|\xi|)^{-\alpha}.
\end{align*}

\textup{(2)}.
 Lemma~\ref{l33} gives that
$\mu\bigl(W_{m,n}(\psi,\bta)\bigr)=\mu(\T^{mn}).$
So it remains to show
$\mu\bigl(W_{m,n}(c\psi,\bta)\bigr)=0$
for every \(0<c<1\). To this end, we claim that
$$\sum_{j=1}^\infty \mu(\mathcal B_j)<\infty.$$
Once the claim is established, by Borel--Cantelli lemma and the obvious inclusion
$$W_{m,n}(c\psi,\bta)\subseteq\limsup_{j\to\infty}\mathcal B_j$$
one has that
$$\mu\bigl(W_{m,n}(c\psi,\bta)\bigr)=0.$$

It remains to prove the claim. Since \(P_j\ge\mathbf1_{\mathcal B_j}\),
by the weak limit of $\{\mu_j\}$ and the stability estimates in \textup{(b)},
\begin{align*}
\mu(\mathcal B_j)\le \int P_j\,\dif\mu&\le\int P_j\,\dif\mu_j+\sum_{\ell=j+1}^\infty\left|\int P_j\,\dif\mu_\ell-\int P_j\,\dif\mu_{\ell-1}\right| \\
&\le\mu_j(\mathcal B_j)+2^{-j-2}+\sum_{\ell=j+1}^\infty2^{-\ell-2}.
\end{align*}

On the other hand, since
$$
\mathcal B_j\subseteq \bigcup_{|\bb|>T_j}\mathcal D_{\eta_j}(\bb),
$$
property \textup{(a)} gives that
$$
\mu_j(\mathcal B_j)\le\prod_{i=1}^{j-1}\|f_i\|_\infty\int f_j(A)\mathbf1_{\bigcup_{|\bb|>T_j}\mathcal D_{\eta_j}(\bb)}(A)\,\dif A<2^{-j-2}.
$$
It follows that 
$$
\mu(\mathcal B_j)\ll2^{-j}, \ \ {\text{and hence}}\ \sum_{j=1}^\infty\mu(\mathcal B_j)<\infty.
$$
 Therefore,
$$
\mu\bigl(\Exact_{m,n}(\psi,\bta)\bigr)=\mu\bigl(\Mat_{m,n}(\T)\bigr),
$$
and \textup{(2)} follows.

Finally, normalizing $\mu$ to a probability measure, 
it gives full mass to $\Exact_{m,n}(\psi,\bta)$ and has Fourier decay
\[
|\widehat\mu(\xi)|\ll(1+|\xi|)^{-\alpha}.
\]
Therefore,
\[
\fdim\Exact_{m,n}(\psi,\bta)\ge2\alpha.
\]
Since $\alpha<n/(1+\tau_\psi)$ can be taken arbitrarily close to $n/(1+\tau_\psi)$, the desired lower bound follows.
\end{proof}

\end{document}